\documentclass[12pt]{amsart}
\usepackage{amsmath}
\usepackage{amssymb}
\usepackage{bm}
\usepackage{graphicx}
\usepackage{verbatim}
\newtheorem{Theorem}{Theorem}
\newtheorem{Lemma}{Lemma}
\newtheorem{Proposition}{Proposition}
\newtheorem{Corollary}{Corollary}

\newtheorem{Conjecture}{Conjecture}

\title[Log-scale Quasimodes that do not Equidistribute]{Logarithmic-scale Quasimodes that do not Equidistribute}
\author{Shimon Brooks}
\thanks{The author was partially supported by NSF grant DMS-1101596 and a Marie Curie Career Integration Grant.}

\begin{document}
\maketitle

{\em Abstract:}  Given any compact hyperbolic surface $M$, and a closed geodesic on $M$, we construct of a sequence of quasimodes on $M$ whose microlocal lifts concentrate positive mass on the geodesic.  Thus, the Quantum Unique Ergodicity (QUE) property does not hold for these quasimodes.  This is analogous to a construction of Faure-Nonnenmacher-De Bi\`evre in the context of quantized cat maps, and lends credence to the suggestion that large multiplicities play a role in the known failure of QUE for certain ``toy models" of quantum chaos.  We moreover conjecture a precise threshold for the order of quasimodes needed for QUE to hold--- the result of the present paper shows that this conjecture, if true, is sharp.

\section{Introduction}\label{intro}

An important problem of ``quantum chaos" is to understand the relationship between a dynamical system and its high-energy quantum mechanical models.  The ``correspondence principle" dictates that quantum mechanics should replicate classical mechanics in  the semiclassical limit $\hbar\to 0$, and this is generally the case when the classical system is completely integrable.  In situations where the dynamics are more complicated, however, the standard analytic techniques are only valid for ``short" times, up to the scale of the Ehrenfest time $C|\log{\hbar}|$, where $C$ is a constant depending on the system.  

Consider the geodesic flow on a compact surface $M=\Gamma\backslash \mathbb{H}$ of constant negative curvature.  It is well known that this dynamical system is ``chaotic"--- eg., mixing (with respect to Liouville measure on $S^*M$), Anosov, etc.  On the other hand, the high-energy spectral data for such surfaces is extremely mysterious (see eg. \cite{SarnakHyp}).  The Quantum Unique Ergodicity (QUE) Conjecture states that high energy eigenfunctions become equidistributed as the eigenvalue tends to $\infty$ --- corresponding to the semiclassical limit $\hbar\to 0$; precisely, that the distributions
$$\mu_\phi : f\in C^\infty(S^*M) \mapsto \langle Op(f)\phi, \phi\rangle$$
converge in the weak-* topology to the Liouville measure on $S^*M$ as the eigenvalue of $\phi$ tends to $\infty$.  Any weak-* limit point of the $\mu_\phi$ is called a {\bf quantum limit}.  

It is known by \cite{Shn, Z2, CdV} that any quantum limit is a (positive) measure invariant under the geodesic flow; if we normalize $||\phi||_2=1$, then they are probability measures.  They further show  that {\em almost all} $\mu_\phi$ become equidistributed, in the following sense:  for any choice of orthonormal basis of $L^2(M)$ consisting of Laplace eigenfunctions, there exists a zero-density exceptional set of basis eigenfunctions such  that
the remaining eigenfunctions satisfy QUE.  This property is known as {\bf Quantum Ergodicity}, and holds in great generality--- it depends only on the ergodicity of the flow.  QUE then asks that there be no exceptional subsequences.

In contrast, there are ``toy models" of quantum chaos that demonstrate non-QUE behavior.  For example, hyperbolic linear maps $A\in SL(2,\mathbb{Z})$ of $\mathbb{T}^2$--- also known as ``cat maps"--- can be quantized, and shown to satisfy Quantum Ergodicity \cite{BouDB}.  Yet not every quantum limit need be Lebesgue; examples are constructed by Faure-Nonnenmacher-De Bi\`evre \cite{FNDB} of quantum limits that are half-Lebesgue and half atomic (i.e., half of its mass is supported on a finite periodic orbit).  

It is suspected \cite{SarnakProgress} that large spectral multiplicities play a role in this phenomenon.  For the cat maps, the quantum propagator has finite order, given by $\text{ord}_N(A)$ where $A\in SL(2,\mathbb{Z})$ is the matrix generating the classical dynamics, $N=(2\pi\hbar)^{-1}\in\mathbb{Z}$ is the inverse Planck's constant, and $\text{ord}_N(A)$ is the order of the matrix $A$ modulo $N$.  It is known that  $\text{ord}_N(A)\geq 2 \log_\lambda{N} + O(1)$, where $\lambda>1$ is the absolute value of the large eigenvalue of $A$.  The examples of \cite{FNDB} rely crucially on minimal periods of the propagator; the sequence of eigenvectors is taken from values of $N$ for which $\text{ord}_N(A)$ saturates this lower bound.  Since the propagator is unitary, the eigenvalues must be $\text{ord}_N(A)$-roots of unity;  the Hilbert space of states from which these eigenvectors are taken is $N$-dimensional, and so such short periods imply massive degeneracies in the spectrum of size
about
$$\frac{N}{2\log_\lambda{N}} \sim \frac{1/2}{\hbar|\log_\lambda{\hbar}|}$$
In fact,  Bourgain shows \cite{Bour} (improving on earlier results of Kurlberg-Rudnick \cite{KR}) that QUE holds whenever $\text{ord}_N(A) \gtrsim_\epsilon N^\epsilon \gtrsim \hbar^{-\epsilon}$.  However, there is still a significant gap between the logarithmic scale $O(1/|\log{\hbar}|)$ of the non-QUE examples and Bourgain's  QUE bounds.

In fact, we conjecture the correct threshold to be just beyond the logarithmic scale of the known counterexamples:
\begin{Conjecture}\label{catconj}
Let $\mu$ be a quantum limit for the cat map $A$, such that the orders of the propagators satisfy
$$\liminf_{N\to\infty} \frac{\text{ord}_N(A)}{\log{N}} \geq C$$
Then the entropy of $\mu$ is bounded below by
$$h(\mu) \geq h_\text{max} - \frac{1}{C} $$

In particular, if 
$$	\liminf_{N\to\infty} \frac{\text{ord}_N(A)}{\log{N}} = \infty	$$
then $\mu$ is Lebesgue measure.
\end{Conjecture}

We now return to the geodesic flow on a compact hyperbolic surface $M$.  We parametrize the spectrum of the Laplacian on $M$ by $\lambda = \frac{1}{4}+r^2$, and here the frequency $r\to\infty$ plays the role of the inverse Planck's constant in the semiclassical limit.  It is conjectured \cite{IwaniecSarnak} that the multiplicities are bounded by $O_\epsilon(r^\epsilon)$, but this problem is well out of the reach of current technology;  in any case we do not have control over different multiplicity scales here as in the case of cat maps.

However, one can artificially introduce ``degeneracies" by considering quasimodes, or approximate eigenfunctions.  Define an $\omega(r)$-{\bf quasimode with approximate parameter} $r$ to be a function $\psi$ satisfying
$$||(\Delta + (\frac{1}{4}+r^2))\psi||_2 \leq r \omega(r)||\psi||_2$$
The  factor of $r$ in our definition comes from the fact that $r$ is essentially the square-root of the Laplace eigenvalue.  Intuitively, we think of $\psi$ as being localized spectrally near the window $[r-\omega(r), r+\omega(r)]$.

The main term in Weyl's Law (see eg. \cite{Berard}) says that the asymptotic density of eigenfunctions near spectral parameter $r$ is proportional to $r$ (with a constant depending on the area of the surface).  Though controlling the error term is a very difficult problem--- indeed, this is precisely the problem of bounding multiplicities--- it is known \cite{Berard} that this approximation is valid for ``large logarithmic windows" $\omega(r)\geq K/\log{r}$, where again $K$ depends on the surface, and is believed to be valid for far smaller windows; eg., $\omega(r) = r^{-\epsilon}$.  In any case, we shall be interested in this paper with windows of size $\omega(r) \gtrsim 1/\log{r}$, where we certainly expect the error in Weyl's Law to be negligible.  By varying the size $\omega(r)$ of the windows, we can analyze its role in QUE phenomena.
\begin{Conjecture}\label{hypconj}
Let $M$ be a compact hyperbolic surface, and $\{\psi_j\}$ a sequence of $o\left(\frac{1}{\log{r}}\right)$-quasimodes.  Then $\{\psi_j\}$ satisfies QUE.
\end{Conjecture}

This conjecture is a slight strengthening of the QUE Conjecture \cite{RS} for the case of Riemann surfaces, which essentially claims (in a quantitative way) that the small spectral multiplicities are responsible for QUE in this setting.  It is also important to observe that the windows $\omega(r) = o(1/\log{r})$ are just beyond what can be analyzed at present:  as remarked by Sarnak \cite{SarnakProgress}, any proof of QUE is expected to address the multiplicity issue (though perhaps indirectly).  Conjecture~\ref{hypconj} suggests that QUE is on par with a multiplicity bound of $o(1/\log{r})$, whereas the current best known bound due to Berard \cite{Berard} is $O(1/\log{r})$.

In this paper, we study ``logarithmic-scale" quasimodes, by which we mean $\left(\frac{\epsilon}{\log{r}}\right)$-quasimodes for some constant $\epsilon>0$.  Precisely, we prove
\begin{Theorem}\label{loc quasi}
Let $M$ be a hyperbolic Riemann surface, and $\xi\subset S^*M$ a closed geodesic.  Then for any $\epsilon>0$, there exists $\delta>0$ and a sequence of $\left(\frac{\epsilon}{\log{r}}\right)$-quasimodes on $M$ whose microlocal lifts do not equidistribute; indeed, they concentrate mass $\geq \delta$ on the geodesic $\xi$.
\end{Theorem}
Thus, if true, Conjecture~\ref{hypconj} is sharp.  This is reminiscent of--- and, indeed, the argument inspired by--- the constructions of \cite{FNDB} for cat maps.

We also remark that the case of arithmetic {\em joint} quasimodes is radically different, as it was shown in \cite{jointQmodes} that any sequence of joint $o(1)$-quasimodes satisfies QUE, even though this includes spaces of much larger dimension $o(r)$ (with arbitrarily slow decay).  This is due to the additional rigid structure imposed by the arithmetic symmetries, which are not respected by the quasimodes  constructed here in Theorem~\ref{loc quasi} (they {\em a fortiori} cannot be $o(1)$-quasimodes of a Hecke operator).

The main part of the argument is constructing $C/\log{r}$-quasimodes that work, for {\em some} large $C$; moving from a large coefficient to a small $\epsilon$ is achieved by means of a (somewhat cheap) argument to dilute out the $C/\log{r}$-quasimode, given in section~\ref{big to small}.  The construction of these quasimodes for a sufficiently large $C$  is the subject of section~\ref{construction}, and is based on a construction from \cite{localized_example}.  In section~\ref{stationary phase estimates} we compute the stationary phase asymptotic which is the analytic heart of the argument, and then apply it in section~\ref{localization main} to prove that our examples are indeed $O\left(\frac{1}{\log{r}}\right)$-quasimodes that localize a positive proportion of their mass on the geodesic $\xi$.

Though the technical approach is slightly altered for simplicity, intuitively the main idea is to begin with a radial $1$-quasimode localized near $p\in\xi$ (more precisely, localized near radial vectors in a neighborhood of $p$), and average over propagation times to improve the quality of the quasimode to order $C/\log{r}$.  To get localization near $\xi$, we then average these spherical quasimodes over a piece of the stable horocycle through $p$, which enhances the quasimodes along $\xi$ through constructive interference, while canceling out the mass away from $\xi$ through destructive interference.  Naturally, then, the main analytic tool for exploiting the interferences will be the method of stationary phase,  which will yield asymptotics for the relevant integrals defining the microlocal lifts of the quasimodes.

{\bf Acknowledgements.}  The question of proving the existence of such quasimodes was first posed to the author by Peter Sarnak.  We also thank Lior Silberman and Elon Lindenstrauss for helpful discussions.

\section{Construction of the Quasimodes}\label{construction}

We begin with a generalization of the discussion in \cite{localized_example}, where we exhibited spherical kernels whose microlocal lifts localize near radial  vectors.  Intuitively, we wish to take a localized kernel  and improve the order of the quasimodes by averaging over propagation times with an operator like
$$\frac{1}{T} \int_0^T e^{-itr} e^{it\sqrt{-\Delta-1/4}}dt$$
which cancels out components of spectral parameter away from $r$.  This is reminiscent of--- and largely inspired by--- the construction in \cite{FNDB} of scarred eigenstates for quantized cat maps.  To get an $O\left(\frac{1}{\log{r}}\right)$-quasimode, we would have to average up to time $T \gtrsim\log{r}$.

From a technical perspective, though, it is easier to proceed in the following manner--- which is essentially just a smoothed version of this averaging procedure.  Recall the Selberg/Harish-Chandra transform for a spherical kernel $k$ on $\mathbb{H}$ (see eg. \cite[Chapter 1.8]{Iwaniec}), giving the eigenvalue $h(s)\phi_s = k\ast \phi_s$ for a Laplace eigenfunction $\phi$ of spectral parameter $s$:
\begin{eqnarray}
h(s) & = & \int_{-\infty}^\infty e^{isu}g(t)dt\nonumber\\
g(t) & = & 2Q\left(\sinh^2\left(\frac{t}{2}\right)\right)\nonumber\\
k(u) & = & -\frac{1}{\pi} \int_u^\infty \frac{dQ(\omega)}{\sqrt{\omega-u}}\label{Selberg/HC}
\end{eqnarray}
The coordinate $u(z,w) = \sinh^2(d(z,w)/2)$ is often more convenient for calculations than the actual distance ($du$ is the radial volume measure on $\mathbb{H}$).  What is most important for our purposes is that whenever $g$ is compactly supported in the interval $[-T, T]$, the kernel $k$ will be supported in the ball of radius $T$ in $\mathbb{H}$.

We can write such a kernel $k$ as a (left-$K$-invariant) function on $\mathbb{H}$, and use Helgason's Fourier inversion \cite{Hel} to write
$$k(z) = \int_{s=0}^\infty \int_{B} e^{(is+\frac{1}{2})<z,b>}\hat{k}(s,b) s\tanh{(\pi s)} dsdb$$
where $b\in B$ runs over the boundary $S^1$ of the disc model for $\mathbb{H}$, and $<z,b>$ represents the (signed) distance to the origin $o$ from the horocycle through the point $z\in\mathbb{H}$ tangent to $b\in B$.  Since each plane wave $e^{(-is+\frac{1}{2})<\cdot,b>}$ is an eigenfunction of spectral parameter $s$, the Fourier transform
\begin{eqnarray*}
\hat{k}(s,b) & = & \int_\mathbb{H} e^{(-is+\frac{1}{2})<z,b>}k(z)dz\\
& = & h(s) e^{(-is+\frac{1}{2})<o,b>}\\
& = & h(s)
\end{eqnarray*}
so that
$$k(z) = \int_{s=0}^\infty \left(  \int_{B} e^{(is+\frac{1}{2})<z,b>}db \right) h(s) s\tanh{(\pi s)} ds$$

It will be more convenient  to write as in \cite{Z2} 
$$e^{(is+\frac{1}{2})<z,b>}db = e^{(is-\frac{1}{2})<z,b>}d\theta =  e^{(is-\frac{1}{2})\varphi(g.k_\theta)}d\theta$$
where $\varphi(g)$ is again the signed distance from the origin to the horocycle through $g\in PSL(2,\mathbb{R})$, and $k_\theta$ parametrizes the $SO(2)$ fibre  $gK$.  Since $\varphi$ is left $K$-invariant and right $N$-invariant, it is sometimes convenient to use $KAN$ coordinates to write
$$g=\begin{pmatrix}  \cos{\theta} & -\sin{\theta}\\ \sin{\theta} & \cos{\theta}\end{pmatrix} \begin{pmatrix} e^{t/2} & 0\\ 0 & e^{-t/2}\end{pmatrix}  \begin{pmatrix}  1 & n\\ 0 & 1\end{pmatrix} = \begin{pmatrix} a& b\\ c& d\end{pmatrix}$$
so that this distance is given by $\varphi(g) = t =\log{(a^2+c^2)}$.

Fix an orthonormal basis $\{\phi_l\}$ of $L^2(M)$ consisting of Laplace eigenfunctions, which we can take to be real-valued for simplicity.  Each eigenfunction generates, under right translations, an irreducible representation $V_{l} = \overline{\{\phi_l(xg^{-1}) : g\in PSL(2,\mathbb{R})\}}$ of $PSL(2,\mathbb{R})$, which together span a dense subspace of $L^2(S^*M)$.  

We distinguish the pairwise orthogonal {\bf weight spaces} $A_{2n}$  in each representation, consisting of those functions satisfying $f(gk_\theta) = e^{ 2i n\theta}f(g)$ for all $k_\theta = \begin{pmatrix} \cos{\theta} & -\sin{\theta}\\ \sin{\theta} & \cos{\theta}\end{pmatrix} \in K$ and $g\in S^*M$.  The weight spaces together span a dense subspace of $V_{l}$.  Each weight space is one-dimensional in $V_{l}$, spanned by $\phi^{(l)}_{2n}$ where 
\begin{eqnarray*}
\phi_0^{(l)} & = & \phi_l \in A_0\\
(ir_l + \frac{1}{2}+n)\phi_{2n+2}^{(l)} & = & E^+\phi_{2n}^{(l)}\\ 
(ir_l + \frac{1}{2} -n)\phi_{2n-2}^{(l)} & = & E^-\phi_{2n}^{(l)}
\end{eqnarray*}
Here $E^+$ and $E^-$ are the {\bf raising} and {\bf lowering operators}, first-order differential operators corresponding to $\begin{pmatrix} 1 & i\\ i & -1\end{pmatrix}\in \mathfrak{sl}(2,\mathbb{R})$ and $\begin{pmatrix} 1 & -i\\ -i & -1 \end{pmatrix}\in \mathfrak{sl}(2,\mathbb{R})$ in the Lie algebra.  
The normalized pseudodifferential operators 
\begin{eqnarray*}
R &  : &  \phi_{2n}^{(l)}\mapsto \phi_{2n+2}^{(l)}\\
R^{-1} & : & \phi_{2n}^{(l)}\mapsto \phi_{2n-2}^{(l)}
\end{eqnarray*}
are unitary and left-invariant, and each $\phi_{2n}$ is a unit vector.
We define the distribution
$$\Phi_\infty^{(l)} = \sum_{n=-\infty}^{\infty} \phi_{2n}^{(l)}$$
and extend this definition by linearity to  $\Psi_\infty = \sum_{n=-\infty}^\infty \psi_{2n} = \sum_{n=-\infty}^\infty R^n\psi$ for quasimodes $\psi$, where each $\psi_{2n} = R^n\psi$.

\subsection{Construction of the Microlocal Lifts}\label{ml lifts construction}
The following can be found in \cite{localized_example}, based on the arguments of \cite{LinHxH}.  We set 
$$I_\psi(f) =\langle Op(f)\psi, \psi\rangle := \langle f \Psi_\infty, \psi\rangle = \lim_{N\to\infty} \left\langle f \sum_{n=-N}^N \psi_{2n}, \psi_0\right\rangle$$
according to the pseudo-differential calculus of \cite{Z2}, which clearly restricts to the measure $|\psi(z)|^2dz$ when applied to $K$-invariant functions $f\in C^\infty(M)$, by orthogonality of the weight spaces.
Note moreover that this limit is purely formal for $K$-finite $f$, and since these $K$-finite functions are dense in the space of smooth functions, we can restrict our attention to these.  

\begin{Lemma}\label{asymp to op}
Let $\psi$ be a linear combination of eigenfunctions with spectral parameter in $[r-1, r+1]$, with $||\psi||_2=1$, and set 
$$\Psi := \sqrt{\frac{3L}{2L^2+1}}\sum_{|n|\leq L} \frac{L - |n|}{L} \psi_{2n}$$  
Then for any smooth $f\in C^\infty(S^*M)$, we have
$$I_\psi(f) = \langle f\Psi, \Psi \rangle +O_{f}(L^{-1} + Lr^{-1})$$
In particular, $\left| I_\psi(f) - \langle f\Psi, \Psi\rangle\right| \to 0$ if we set $L=\lfloor r^{\alpha}\rfloor $ for some $0<\alpha<1$.
\end{Lemma}

{\em Proof:}  See \cite{LinHxH, jointQmodes, localized_example}.

Note that the prefactor $\sqrt{\frac{3L}{2L^2+1}}\sim \sqrt{\frac{3}{2}}r^{-1/4}$ is simply an $L^2$-normalization of the Fej\'er coefficients $\frac{L - |n|}{L} $.

\begin{Corollary}
Let $\{\psi_j\}$ be a sequence of $o(1)$-quasimodes with approximate spectral parameters $r_j$, normalized in $L^2(M)$, and $L=r_j^{-\alpha}$ for some fixed $0<\alpha<1$ as above.  Then for all $f\in C^\infty(S^*M)$ we have
$$\left| I_{\psi_j}(f) -  \langle f\Psi_j, \Psi_j\rangle\right| \to 0$$
\end{Corollary}

{\em Proof:}  As in the proof of Lemma~\ref{asymp to op} from \cite{LinHxH, localized_example}, we may take $f\in\sum_{|n|\leq N_0}A_{2n}$ to be $K$-finite, since these span a dense subspace of $C^\infty(S^*M)$.  The Lemma shows that the statement holds for the projection of $\psi_j$ to the space spanned by eigenfunctions of spectral parameter in $[r_j-1,r_j+1]$; it is sufficient to check that the contribution of other spectral components is negligible.  We decompose $\psi = \psi_\text{in}+\psi_{\text{out}}$, and the corresponding $\Psi = \Psi_\text{in} + \Psi_\text{out}$, where the ``in" component consists of all spectral components inside the interval $[r-1, r+1]$, and the ``out" component is the orthogonal complement consisting of spectral components outside this range.  

We have
\begin{eqnarray*}
\langle f \Psi, \Psi\rangle - \langle f \Psi_\text{in}, \Psi_\text{in}\rangle 
& \leq & \Big| \langle f\Psi, \Psi_{\text{out}} \rangle\Big| +  \Big| \langle f\Psi_{\text{out}}, \Psi_{\text{in}}\rangle\Big| \\
& \leq & ||f||_\infty \cdot \left(\left\|  \Psi\right\|_2 \cdot \left\| \Psi_\text{out} \right\|_2  +  \left\| \Psi_{\text{out}}\right\|_2\cdot \left\| \Psi_{\text{in}}\right\|_2 \right)\\
& \lesssim_f & ||\Psi_\text{out}||_2 = ||\psi_\text{out}||_2\\
I_\psi(f) - I_{\psi_\text{in}}(f) & \leq & \left| \left\langle f \sum_{|n|\leq N_0} \psi_{2n}, \psi_{\text{out}}\right\rangle\right| + \left| \left\langle f\sum_{|n|\leq N_0} \left(\psi_{\text{out}}\right)_{2n}, \psi_{\text{in}}\right\rangle\right| \\
& \leq & (2N_0+1) \cdot ||f||_\infty \cdot \left(||\psi||_2 \cdot  ||\psi_\text{out}||_2 + ||\psi_{\text{out}}||_2 \cdot ||\psi_{\text{in}}||_2\right) \\
& \lesssim_f & ||\psi_\text{out}||_2
\end{eqnarray*}
recalling that each $||\psi_{2n}||_2=||\psi||_2$ by unitarity of the raising and lowering operators $R$ and $R^{-1}$, and so it remains to show that $||\psi_\text{out}||_2 = o(1)$.  But note that for any spectral component $\phi_l$ with spectral parameter outside $[r-1, r+1]$, the defect $||(\Delta - (\frac{1}{4}+r^2))\phi_l||_2 \gtrsim r||\phi_l||_2$, and thus
$$||(\Delta - (\frac{1}{4}+r^2))\psi||_2 \gtrsim ||\psi_\text{out}||_2\cdot r$$
so that the $o(1)$-quasimode condition necessitates $||\psi_\text{out}||_2=o(1)$, and we get the desired asymptotic.  $\Box$

Thus, for any given sequence $\{\psi_j\}$ of $o(1)$-quasimodes (and in particular, for sequences of $O(1/\log{r})$-quasimodes), we have constructed a sequence $\{\Psi_{j}\}_{j=1}^\infty$ such that the microlocal lifts $|\Psi_{j}|^2dVol$ are asymptotically equivalent to the distributions $I_{\psi_j}$.  It is these measures that we wish to study.

\subsection{The kernels $\kappa_{\xi,j}$}\label{kernel}

Propagating spherical kernels over long times will equidistribute on $S^*M$, but by averaging over a stable neighborhood of the geodesic $\xi$, we will obtain quasimodes that localize a positive proportion of their mass on $\xi$.  The intuitive picture is that averaging along a stable horocycle causes constructive interference at vectors pointing along $\xi$, enhancing the mass near the geodesic, while causing destructive interference away from $\xi$.  Naturally, the basis for our analysis in the coming sections will be the method of stationary phase.

Thus, we pick a point $p\in\xi$ on our geodesic, and choose a covering $\mathbb{H}\to M$ such that the origin $i\in\mathbb{H}$ is mapped to $p$, and the imaginary axis $x=0$ is mapped onto the geodesic $\xi$.  Let $\{r_j\}\subset l(\xi)\cdot 2\pi\mathbb{Z}$ be a sequence of resonant frequencies for the length $l(\xi)$ of the geodesic $\xi$; this condition on the approximate spectral parameters will ensure that our quasimodes will not self-interfere as they wrap around $\xi$, and will be used in the calculations of Proposition~\ref{Uj big}. 

We then take the spherical kernel given, via (\ref{Selberg/HC}), by the Fourier pair 
\begin{eqnarray*}
\tilde{h}(s) & = &  \frac{\cosh{\frac{s}{2K_j}}\cosh{\frac{{r_j}}{2K_j}}}{\cosh{\frac{s}{K_j}}+\cosh{\frac{r_j}{K_j}}}\\
\tilde{g}(\xi) & = &  \frac{1}{2}K_j  \frac{\cos(\xi r_j)}{\cosh\left(\pi K_j\xi\right)}
\end{eqnarray*}
for $K_j : = \frac{C}{2\log{r_j}}$.
We then choose a smooth, positive, even cutoff function $\chi$ supported in $[-1,1]$, which we may assume is identically $1$ on $[-\frac{1}{2}, \frac{1}{2}]$, and replace 
\begin{eqnarray*}
g(\xi) & = & \tilde{g}(\xi)\chi\left(	\frac{C}{\log{r_j}} \xi	\right)\\
h(s) & = & \tilde{h}(s) \ast \hat\chi\left(	\frac{\log{r_j}}{C} s	\right)
\end{eqnarray*}
so that $g$ is supported in the interval of radius $C^{-1}\log{r_j}$.  We denote the resulting kernel by $k_j$,
and its microlocal lift--- as constructed in Lemma~\ref{asymp to op}--- by $\kappa_j$; they are supported in the ball of radius $C^{-1}\log{r_j}$ in $\mathbb{H}$ and $S^*\mathbb{H}$ respectively.   Note as in \cite{localized_example} that the distribution
$$	\sum_{n=-\infty}^\infty R^nk_j (g) = \int_s e^{(is-\frac{1}{2})\varphi(g.k_\theta)}d\theta h(s) s \tanh(\pi s) ds	$$
by \cite{Z2}, with $k_\theta = \begin{pmatrix} \cos\theta & -\sin\theta\\ \sin\theta & \cos\theta \end{pmatrix}$, and so we may write $\kappa_j$ as a (right-)convolution of $\sum_{n=-\infty}^\infty R^n k_j$ with the $L^2$-normalized Fej\'er kernel of order $L$ on the $SO(2)$-fibre $K$:
$$\kappa_j = \int_s \int_{\theta=0}^\pi e^{(is-\frac{1}{2})\varphi(g.k_\theta)}F_L(2\theta)d\theta h(s) s \tanh(\pi s) ds$$
where $F_L(2\theta):= \sqrt{\frac{3L}{2\pi(2L^2+1)}}\sum_{|n|\leq L}\frac{L-|n|}{L} \cos{(2n\theta)}$.
We will take $L=\lfloor r_j^{10/C}\rfloor$ as a convenient\footnote{It should be emphasized that here and throughout we make no attempt to obtain the sharpest possible bounds, and thus allow ourselves to be wasteful with powers of $r_j^{1/C}$, at the expense of taking a larger constant  $C$ than necessary.} choice for $L$.

We note for later use the estimate
\begin{equation}\label{total spec}
\int_s s^{-1}h(s) \tanh{\pi s } ds  \lesssim  \frac{1}{r_j \log{r_j}}\lesssim r_j^{-1}
\end{equation}
and remark, as in \cite{localized_example}, that since $\chi\in C^\infty$ is smooth, combined with the rapid decay of $\tilde{h}$ away from $r_j$, the function $h(s)$ is also rapidly decaying away from $r_j$; the estimate $h(s) \lesssim |r_j-s|^{-3}$ will suffice here.

We now define 
$$k_{\xi,j}(z) = \int_n k_j\left(\begin{pmatrix}1& n\\ & 1\end{pmatrix}.z\right) \chi(n)dn$$
 to be the convolution of $k_j$ with the smooth cutoff $\chi$, along the stable subgroup $N$.  Since this convolution acts in $k_j$ on the left, it commutes with the left-invariant operators $R$ and $R^{-1}$ defining our microlocal lift, and so $\kappa_{\xi, j}$ is given by
\begin{eqnarray*}
\kappa_{\xi, j} (g) & = & \int_n \chi(n)\kappa(n.g) dn\\
& = & \int_s\left(\int_{n} \int_{\alpha\in\mathbb{T}} \chi(n) F_L(2\alpha)  e^{(is-\frac{1}{2})\varphi(n.g.k_\alpha)} d\alpha dn\right)h(s)s\tanh(\pi s) ds
\end{eqnarray*}
We also define 
\begin{eqnarray*}
\bar{k}_{\xi,j}(z) & = & \sum_{\gamma\in\Gamma}k_{\xi,j}(\gamma.z)\\
\bar{\kappa}_{\xi,j} (g) & = & \sum_{\gamma\in\Gamma}\kappa_{\xi,j}(\gamma.g)
\end{eqnarray*}
to be the respective projections to $M$ and $S^*M$.  Note that since $k_{\xi,j}$ and $\kappa_{\xi,j}$ have compact support, the sums are finite; and that as usual $||\bar{k}_{\xi,j}||_{L^2(M)} = ||\bar{\kappa}_{\xi,j}||_{L^2(S^*M)}$.  Naturally these projections do not commute with the averaging over the stable horocycle that defines $k_{\xi,j}$ and $\kappa_{\xi,j}$, but it will still be easier for us to approximate $\kappa_{\xi,j}$ pointwise, and then estimate the sum over $\Gamma$ to project back to $S^*M$; in reality, the vast majority of the values of $\kappa_{\xi,j}$ in the sum over $\Gamma$ will be negligible.

We will show that:
\begin{itemize}
\item{For any neighborhood $U$ of the geodesic $\xi$ in $S^*M$, we have $\int_U |\bar{\kappa}_{\xi,j}|^2 dVol \gtrsim 1/\log{r_j}$.}
\item{The full $L^2$-norm satisfies $||\bar{\kappa}_{\xi, j}||_2^2 \lesssim 1/\log{r_j}$.}
\item{The function $\bar{k}_{\xi,j}$ is an $O\left(\frac{1}{\log{r_j}}\right)$-quasimode on $M$.}
\end{itemize}
Here and throughout, we allow all implied constants to depend on $C$.
Together these prove Theorem~\ref{loc quasi} for the case of a large constant; in the final section~\ref{big to small} we show how to bootstrap the argument to $\left(\frac{\epsilon}{\log{r_j}}\right)$-quasimodes, for arbitrarily small constants $\epsilon>0$.

\section{A Stationary Phase Estimate}\label{stationary phase estimates}

We now perform the main stationary phase approximation that is the crux of our analysis.  We write $\kappa^{(s)}_{\xi,j}$ for the $s$-spectral component of $\kappa_{\xi,j}$; that is to say
\begin{eqnarray*}
\kappa_{\xi,j} & = & \int_s \kappa^{(s)}_{\xi,j} h(s) s \tanh(\pi s) ds\\
\kappa^{(s)}_{\xi, j} (x,t,\theta) & = & \int_{n} \int_{\alpha\in\mathbb{T}} \chi(x-n) F_L(2|\theta-\alpha|)  e^{(is-\frac{1}{2})\varphi(n,t,\alpha)} d\alpha dn
\end{eqnarray*}
where we denote
\begin{eqnarray}
\varphi(n,t,\alpha) & = & \varphi\left(	\begin{pmatrix} 1&n\\ & 1\end{pmatrix} \begin{pmatrix} e^{t/2}& \\ & e^{-t/2}\end{pmatrix} \begin{pmatrix}	\cos\alpha & -\sin\alpha\\ \sin\alpha & \cos\alpha \end{pmatrix}	\right)\nonumber\\
& = & \varphi \begin{pmatrix}	e^{t/2}\cos\alpha + ne^{-t/2}\sin\alpha & \ast\\ e^{-t/2}\sin\alpha & \ast	\end{pmatrix}\nonumber\\
& = & \log(e^t\cos^2\alpha + n\sin(2\alpha)+e^{-t}(n^2+1)\sin^2\alpha)\label{NAK varphi}
\end{eqnarray}
in $NAK$ coordinates.  We will perform the stationary phase analysis in the integrals over $n$ and $\alpha$, estimating $\kappa^{(s)}_{\xi, j}$ pointwise, and then deal with the integral over the spectrum afterwards.

\begin{Proposition}\label{stationary phase estimate}
Let $(x,t,\theta)\in S^*\mathbb{H}$ with $|t|\leq C^{-1}\log{r_j}$, and set $L=\lfloor r_j^{10/C}\rfloor$.    Then
\begin{eqnarray*}
\kappa_{\xi, j}^{(s)}(x,t,\theta) & = & \pi s^{-1} \chi(x) e^{t/2} \left[  e^{ist}F_L(2\theta) + e^{-ist}F_L(\pi-2\theta)	\right] + O(r_j^{100/C}s^{-2})
\end{eqnarray*}
\end{Proposition}
{\em Remark:}  Along the geodesic, the $e^{t/2}$ term represents the growth due to constructive interference in the stable direction.  This exponential term will play a central role in the mass concentration of Proposition~\ref{Uj big}.

{\em Proof:}  The proof is an application of the method of stationary phase (see eg. \cite{Hormander1, EvansZworski}).  Note that $\kappa_{\xi,j}^{(s)}$ is defined pointwise via integrals over $N$ and $K$; the integral over $N$ takes place on the left, while the integral over $K$ takes place on the right, so that these translations commute.  The sole critical points of the phase function $\varphi$ are at the origin $(n,\alpha)=(0,0)$ and at $(n,\alpha)=(0,\pi/2)$; to see this it is actually easiest to work in $KAN$-coordinates, to take advantage of the fact that $\varphi$ is left-$K$ and right-$N$ invariant.  Let $g = \begin{pmatrix} \cos\theta & -\sin\theta\\ \sin\theta & \cos\theta\end{pmatrix} \begin{pmatrix} e^{t/2}& \\ & e^{-t/2} \end{pmatrix} \begin{pmatrix} 1 & u \\ & 1\end{pmatrix} \in PSL(2,\mathbb{R})$; we wish to calculate the critical points of $\varphi$ for translation along $N$ on the left, and $K$ on the right.  Thus we compute
\begin{eqnarray*}
\left. \frac{\partial}{\partial n} \varphi(n.g)\right|_{n=0} & = & 
\left. \frac{\partial}{\partial n}\varphi\left(\begin{pmatrix} 1 & n\\ 	& 1\end{pmatrix}\begin{pmatrix} \cos\theta & -\sin\theta\\ \sin\theta & \cos\theta\end{pmatrix} \begin{pmatrix} e^{t/2}& \\ & e^{-t/2} \end{pmatrix} \begin{pmatrix} 1 & u \\ & 1\end{pmatrix}	\right) \right|_{n=0}\\
& = & \left. \frac{\partial}{\partial n}\varphi\left(\begin{pmatrix} 1 & n\\ 	& 1\end{pmatrix}\begin{pmatrix} \cos\theta & -\sin\theta\\ \sin\theta & \cos\theta\end{pmatrix} \begin{pmatrix} e^{t/2}& \\ & e^{-t/2} \end{pmatrix} \right)\right|_{n=0}\\
& = & \left. \frac{\partial}{\partial n} \varphi \begin{pmatrix}	e^{t/2}(\cos\theta + n\sin\theta) & \ast\\ e^{t/2}\sin\theta & \ast 	\end{pmatrix}   \right|_{n=0} \\
& = & \frac{\partial}{\partial n} \log\Big(e^t(\cos\theta + n\sin\theta)^2 + e^t\sin^2\theta \Big)\Big|_{n=0}\\
& = & \left.  \frac{\partial}{\partial n} (t + \log(1+n\sin2\theta + n^2\sin^2\theta))\right|_{n=0}\\
& = & \left. \frac{\sin2\theta + 2n\sin^2\theta}{1+n\sin(2\theta) + n^2\sin^2\theta}\right|_{n=0}\\
& = & \sin(2\theta)
\end{eqnarray*}
Thus the critical points must satisfy $\sin(2\theta)=0$, which means $\theta=0$ or $\theta=\pi/2$.
Similarly we compute
\begin{eqnarray*}
\left. \frac{\partial}{\partial \alpha} \varphi(g.k_\alpha)\right|_{\alpha=0} & = & 
\left. \frac{\partial}{\partial \alpha}\varphi\left(\begin{pmatrix} \cos\theta & -\sin\theta\\ \sin\theta & \cos\theta\end{pmatrix} \begin{pmatrix} e^{t/2}& \\ & e^{-t/2} \end{pmatrix} \begin{pmatrix} 1 & u \\ & 1\end{pmatrix} \begin{pmatrix}\cos\alpha & -\sin\alpha\\ \sin\alpha & \cos\alpha \end{pmatrix}	\right) \right|_{\alpha=0}\\
& = & \left. \frac{\partial}{\partial \alpha}\varphi\left( \begin{pmatrix} e^{t/2}& \\ & e^{-t/2} \end{pmatrix}\begin{pmatrix} 1& u\\ & 1 \end{pmatrix}\begin{pmatrix}\cos\alpha & -\sin\alpha\\ \sin\alpha & \cos\alpha \end{pmatrix} \right)\right|_{\alpha=0}\\
& = & \left.		\frac{\partial}{\partial\alpha} \varphi\begin{pmatrix} e^{t/2}(\cos\alpha + u\sin\alpha) & \ast\\ e^{-t/2}\sin\alpha & \ast 	\end{pmatrix}	\right|_{\alpha=0}\\
& = & \left.  \frac{\partial}{\partial \alpha} \log\Big(e^t(\cos\alpha + u\sin\alpha)^2 + e^{-t}\sin^2\alpha \Big)\right|_{\alpha=0}\\
& = & 2u
\end{eqnarray*}
and so the critical points of $\varphi$ are $(u,\theta)=(0,0)$ and $(u,\theta)=(0,\pi/2)$ (in $KAN$-coordinates); note however that these are simply the points based at $ie^t \in \mathbb{H}$, pointing to $\infty$ and to $0$--- these points also have coordinates $(0,t,0)$ and $(0,t,\pi/2)$ in $NAK$-coordinates.

\vspace{.1in}

Now recall that
$$\kappa^{(s)}_{\xi, j}(x,y,\theta) = \int_{n} \int_{\alpha\in\mathbb{T}} \chi(x-n) F_L(2|\theta-\alpha|)  e^{(is-\frac{1}{2})\varphi(n,t,\alpha)} d\alpha dn$$
We wish to apply the method of stationary phase to the integral over $\alpha, n$.  
For this we will also need to compute the determinant $|\partial^2\varphi|$ of the Hessian  at the critical points $(0,t,0)$ and $(0,t, \pi/2)$.

First we consider the point $(0,t,0)$, based on the imaginary axis and pointing towards $\infty$.  Since $\varphi$ is left-$N$-invariant for vectors pointing towards $\infty$, the second derivative $\left. \frac{\partial^2\varphi}{\partial n^2}\right|_{\alpha=0} = 0$, and thus the second derivative along $\alpha$ is irrelevant; indeed, the determinant $|\partial^2\varphi|$ of the Hessian is simply $\left|\frac{\partial^2\varphi}{\partial n\partial\alpha}	\right|^2$.  

This derivative is readily  computed using (\ref{NAK varphi}); taking derivatives first with respect to $n$ (at $n=0$) and then with respect to $\alpha$ gives
\begin{eqnarray*}
\left.\frac{\partial^2\varphi}{\partial\alpha\partial n}\right|_{(0,t,0)} & = & \left. \frac{\partial}{\partial\alpha}\right|_{\alpha=0}  \frac{\sin(2\alpha)}{e^t\cos^2\alpha+e^{-t}\sin^2\alpha}\\
& = & \left. \frac{2\cos(2\alpha)}{e^t\cos^2\alpha+e^{-t}\sin^2\alpha} - \frac{\sin^2(2\alpha)(e^{-t}-e^t)}{(e^t\cos^2\alpha+e^{-t}\sin^2\alpha)^2}\right|_{\alpha=0}\\
& = & 2e^{-t}
\end{eqnarray*}
and hence
\begin{equation}\label{up d2}
\left|\partial^2\varphi\right|^{-1/2} \Big|_{(0,t,0)} = \frac{1}{2}e^{t}
\end{equation}

At the point $(0,t, \pi/2)$, on the other hand, we find that
\begin{eqnarray*}
\left. \frac{\partial^2\varphi}{\partial\alpha\partial n}\right|_{(0,t,\pi/2)}
& = & \left. \frac{2\cos(2\alpha)}{e^t\cos^2\alpha+e^{-t}\sin^2\alpha} - \frac{\sin^2(2\alpha)(e^{-t}-e^t)}{(e^t\cos^2\alpha+e^{-t}\sin^2\alpha)^2}\right|_{\alpha=\frac{\pi}{2}}\\
& = & -2e^{t}
\end{eqnarray*}
Here the derivatives along $N$ are no longer identically zero, and so we must also calculate 
\begin{eqnarray*}
\left. \frac{\partial^2\varphi}{\partial\alpha^2}\right|_{(0,t,\pi/2)} & = & \left. \frac{\partial^2}{\partial\alpha^2} \log(e^t\cos^2\alpha + e^{-t}\sin^2\alpha)\right|_{\alpha=\pi/2}\\
& = &\left. \frac{\partial}{\partial\alpha} \left(\frac{-2\sinh(t)\sin(2\alpha)}{e^t\cos^2\alpha + e^{-t}\sin^2\alpha}\right)\right|_{\alpha=\pi/2}\\
& = & 4\sinh(t) e^{t} = 2e^{2t}-2\\
\left. \frac{\partial^2\varphi}{\partial n^2}\right|_{(0,t,\pi/2)} & = &  \left. \frac{\partial^2}{\partial n^2} \log(e^{-t}(n^2+1))\right|_{n=0}\\
& = & \left. \frac{\partial}{\partial n} \left(\frac{2n}{n^2 +1}\right)\right|_{n=0}\\
& = & 2
\end{eqnarray*}
and thus 
\begin{eqnarray}
\left|\partial^2\varphi\right|^{-1/2} \Big|_{(0,t,\pi/2)} & = & |4e^{2t}-4 -4e^{2t}|^{-1/2}\nonumber\\
& = & \frac{1}{2}\label{down d2}
\end{eqnarray}

Thus, recalling that $\varphi(0,t,0)=t$ and $\varphi(0,t,\pi/2)=-t$, we find by the method of stationary phase (eg., \cite[Theorem 7.7.5]{Hormander1})
\begin{eqnarray*}
\kappa^{(s)}_{\xi, j} (x,t,\theta) & = & \int_{n} \int_{\alpha\in\mathbb{T}} \chi(x-n) F_L(2|\theta-\alpha|)  e^{(is-\frac{1}{2})\varphi(n,t,\alpha)} d\alpha dn \\
& = & \frac{2\pi s^{-1}}{\sqrt{|\partial^2\varphi(0,t,0)|}} e^{(is-\frac{1}{2})t}\chi(x)F_L(2\theta) \\
 & & \quad + \quad \frac{2\pi s^{-1}}{\sqrt{|\partial^2\varphi(0,t,\frac{\pi}{2})|}} e^{(is-\frac{1}{2})(-t)}\chi(x)F_L(\pi - 2\theta) + O_{L,t}(s^{-2})\\
& = & \frac{2\pi s^{-1}}{2e^{-t}} e^{(is-\frac{1}{2})t}\chi(x)F_L(2\theta) + \frac{2\pi s^{-1}}{2} e^{(is-\frac{1}{2})(-t)}\chi(x)F_L(\pi - 2\theta) + O_{L,t}(s^{-2})\\
 & = & \frac{\pi e^{t/2}\chi(x)}{s} \Big(e^{ist}F_L(2\theta) + e^{-ist}F_L(\pi - 2\theta)\Big) + O_{L,t}(s^{-2})
\end{eqnarray*}
The error term depends on $2$ derivatives of $F_L$--- which are bounded by $O(L^{3})$---  and of $\chi$, which are fixed.   The dependence on $t$ is bounded by a fixed number of derivatives of $\varphi(\cdot, t, \cdot)$, and powers of $|\partial\varphi(\cdot, t, \cdot)|^{-1}$ away from the critical points, along with the factor $e^{t/2}$; each of which are controlled by powers of $e^{|t|}$.  Thus since $L\lesssim r_j^{10/C}$ and $e^{|t|}\leq r_j^{1/C}$, we can reduce the error term to $O(s^{-2})$ times a fixed power of $r_j^{1/C}$; a crude bookkeeping of the various remainders shows that the estimate $O(s^{-2}r_j^{100/C})$ is more than sufficient.  $\Box$

\vspace{.2in}

Now we continue by performing the integration over the spectrum:
\begin{Proposition}\label{spec integral}
We have
\begin{eqnarray*}
\kappa_{\xi,j}(x,t,\theta) & = & \pi\chi(x)\Big[F_L(2\theta) + F_L(\pi - 2\theta)\Big]e^{t/2}g(t) + O(r_j^{100/C-1})\\
\end{eqnarray*}
\end{Proposition}

{\em Proof:}
Proposition~\ref{stationary phase estimate} shows that
$$\kappa_{\xi, j} (x,t,\theta)= \pi\chi(x)e^{t/2}\int_s h(s)\tanh(\pi s) \left(e^{ist}F_L(2\theta)  + e^{-ist}F_L(\pi-2\theta) \right)ds + O(r_j^{100/C -1})$$
since by (\ref{total spec}) the integral of the error term $O(r_j^{100/C} s^{-2})$ certainly satisfies
\begin{eqnarray*}
r_j^{100/C} \int_s s^{-1}h(s) \tanh(\pi s)ds & \lesssim & r_j^{100/C-1}
\end{eqnarray*}
 Moreover, since $1 > \tanh(\pi s) > 1 - 2e^{-2\pi s}$, and for $|s|\leq \frac{1}{2}r_j$ we have $|h(s)| \lesssim |r_j-s|^{-3}$ by the rapid decay of $h$, we estimate
\begin{eqnarray*}
\left|g(t) - \int_{s}e^{ist}h(s)\tanh(\pi s) ds \right| & = & \left|	\int_s e^{ist}h(s)[1-\tanh(\pi s)]	ds\right|\\
& \lesssim & \int_s h(s)[1-\tanh(\pi s)] ds\\
& \lesssim & \int_{|s|\leq \frac{r_j}{2}} h(s)ds + \int_{|s|\geq \frac{r_j}{2}} h(s) \cdot (1-\tanh(\pi s)) ds\\
& \lesssim &  \int_{|s|\leq \frac{r_j}{2}} |r_j-s|^{-3} ds + \int_{|s|\geq \frac{r_j}{2}} h(s) e^{-2\pi s} ds\\
& \lesssim & r_j^{-2} + e^{-\pi r_j} \int_{|s|\geq \frac{r_j}{2}} h(s)\\
& \lesssim & r_j^{-2}
\end{eqnarray*}
which is readily absorbed into the error term $O(r_j^{100/C -1})$.

Similarly for the $e^{-ist}F_L(\pi-2\theta)$ term, noting that $g$ is even so that 
$$\int_s e^{-ist}h(s)ds = g(-t) = g(t)$$
and we are done.  $\Box$

\begin{Corollary}\label{kappa small}
Assume $C$ is sufficiently large (eg. $C>200$ is sufficient).  Then for $|t|\leq 2C^{-1}\log{r_j}$ and $|2\theta|, |\pi-2\theta| \geq r_j^{-5/C}$, we have
$$|\kappa_{\xi,j}(x,t,\theta)|\lesssim r_j^{-3/C}$$
\end{Corollary}

{\em Proof:}  Plug into Proposition~\ref{spec integral} the estimate 
$$F_L(2\theta) \lesssim L^{-3/2}\frac{\sin^2(n\theta)}{\sin^2\theta}  \lesssim  r_j^{-5/C}$$
whenever $|2\theta|\geq r_j^{-5/C}$ and $L\gtrsim r_j^{10/C}$, and similarly for $F_L(\pi-2\theta)$.   Since $\chi$ and $g(t)$ are 
uniformly bounded and $e^{t/2}\leq r_j^{1/C}$, the main terms of the $\kappa_{\xi,j}$ asymptotic from Proposition~\ref{spec integral} are $O(r_j^{-4/C})$.  As long as $C$ is large, the error term $O(r_j^{100/C-1})$ is also small enough and the Corollary holds.  $\Box$

\section{Localization and the Proof of Theorem~\ref{loc quasi}}\label{localization main}

Recall we have a closed geodesic $\xi\subset M$, and a point $p\in\xi$.  We lift $M$ to the upper-half plane in such  a way that $p$ is mapped to the origin, and $\xi$ is mapped to the imaginary axis $\{iy : y>0\}$.  We will often abuse notation and also use $\xi$ to refer to the geodesic path in $S^*M$.

We consider the ``collars" $U_j\subset PSL(2\mathbb{R})\cong NAK$ defined by
\begin{eqnarray*}
U_j^{\text{up}} & = & \{(x,t,\theta) : |x| \leq 1, \quad |t|\leq C^{-1}\log{r_j}, \quad |2\theta|\leq r_j^{-5/C}\}\\
U_j^{\text{down}} & = & \{(x,t, \theta) : |x| \leq 1, \quad |t|\leq C^{-1}\log{r_j}, \quad |\pi - 2\theta|\leq r_j^{-5/C}\}\\
U_j &=& U_j^{\text{up}} \cup U_j^{\text{down}}
\end{eqnarray*}
for $C$ sufficiently large.

We show that the projection of $U_j$ to $S^*M$ is one-to-one, with the obvious exceptions of the periodicity of the geodesic $\xi$:
\begin{Lemma}\label{almost 1-1}
Let $U_j$ as above.  Then
$$(U_j \cdot U_j^{-1}) \cap \Gamma \subset A$$
\end{Lemma} 
In other words, if $g_1,g_2\in U_j$ project to the same point $\bar{g}\in S^*M$, then $g_1\in\Gamma.g_2$ and thus $g_1g_2^{-1}\in\Gamma$; the lemma then guarantees that $g_1 = a_{k l(\xi)}.g_2$, where $l(\xi)$ is the length of the geodesic $\xi$, and the element of $\Gamma$ in question is simply wrapping around the geodesic $\xi$ exactly $k$ times.

{\em Proof:}  We give the argument for $U_j^{\text{up}}$, the argument for $U_j^{\text{down}}$ is completely analogous (as is the argument that the projections of $U_j^{\text{up}}$ and $U_j^{\text{down}}$ are disjoint).

First, we note that $U_j^{\text{up}} \subset B^-(2)A B^+(r_j^{-4/C})$, where 
$B^-(2)$ is a $2$-ball along the stable direction $N$, while $B^+(r_j^{-4/C})$ is a small $r_j^{-4/C}$ ball in the unstable direction $\bar{N}$, and $A$ is the diagonal subgroup.  Since in  $U_j$ the diagonal component is at most $C^{-1}\log{r_j}$ and the $K$-component is very small $\leq r_j^{-5/C}$, this means that up to minor adjustments to the constants, the $K$-component can be replaced by  $NA\bar{N}$ coordinates that are simpler to manipulate, since we intend to conjugate by large elements of $A$.

Thus, suppose we have
$$n(u_-) a_t \bar{n}(u_+) \in \Gamma n(v_-)a_s \bar{n}(v_+)$$
where the stable components $|u_-|, |v_-| < 2$, the unstable components $|u_+|, |v_+| <r^{-4/C}$, and $|t|, |s| \leq 2C^{-1}\log{r_j}$.  We then conjugate by $a_{k\cdot l(\xi)}\in \Gamma$, for a large integer $k$ of size $\frac{3C^{-1}}{l(\xi)}\log{r_j}$.  This implies
$$a_t n(e^{-kl(\xi)+t}u_-) \bar{n}(e^{kl(\xi)}u_+) \in \Gamma a_s n(e^{-kl(\xi)+s}v_-) \bar{n}(e^{kl(\xi)}v_+)$$
and moreover since $a_{l(\xi)}\in \Gamma$ we can write
$$a_\tau n(e^{-kl(\xi)+t}u_-) \bar{n}(e^{kl(\xi)}u_+) \in \Gamma a_\sigma n(e^{-kl(\xi)+s}v_-) \bar{n}(e^{kl(\xi)}v_+)$$
for $|\tau|, |\sigma|\leq l(\xi)= O(1)$.  But note that 
\begin{eqnarray*}
\left|e^{-kl(\xi)+t}u_- \right|, \left|e^{-kl(\xi)+s}v_- \right| & \lesssim & r_j^{-3/C} r_j^{2/C} \lesssim r_j^{-1/C}\\
\left|e^{kl(\xi)}u_+\right| , \left|e^{kl(\xi)}v_+\right| & \lesssim & r_j^{3/C} r_j^{-4/C} \lesssim r_j^{-1/C}
\end{eqnarray*}
are all small, and so this element of $\Gamma$ is $r_j^{-1/C}$-close (on the right) to the geodesic $\xi$.  This implies that it is in fact on the geodesic itself, once $r_j$ is large, since $\xi$ has a (right-)neighborhood in $S^*M$ that does not intersect any  elements of $\Gamma$ other than those in $A$ (corresponding to the periodicity of $\xi$).  $\Box$

\begin{Proposition}\label{Uj big}
Let $U\supset \xi$ be an open set in $S^*M$ containing the geodesic $\xi$.  Then for $r_j$ sufficiently large,
$$||\bar\kappa_{\xi,j}||_{L^2({U})}^2 \gtrsim 1/\log{r_j}$$
\end{Proposition}

{\em Proof:}  
The neighborhood $U$ contains a subset $V\subset U$ of the form
$$V = \left\{ \xi.\begin{pmatrix} 1 & y\\ & 1\end{pmatrix} \begin{pmatrix} \cos\theta & -\sin\theta\\ \sin\theta & \cos\theta \end{pmatrix} : |y|\leq \eta, \quad |\theta|\leq r_j^{-5/C} \right\}$$
for some $\eta>0$, where we set $\theta=0$ to be a direction along the geodesic contained in $U$.  In this parameterization $V$ is contained in the projection $\bar{U}_j^{\text{up}}$ of $U_j^{\text{up}}$ to $S^*M$, and disjoint from the projection $\bar{U}_j^{\text{down}}$ of $U_j^{\text{down}}$.  Also note
 that since the contributions of $\kappa_{\xi,j}$ from points outside $U_j^{\text{up}}$ and $U_j^{\text{down}}$ are each bounded by $r_j^{-3/C}$ by Corollary~\ref{kappa small}, and since there are $O(r_j^{2/C})$ such points in the support of $\kappa_{\xi,j}$ --- contained in the ball of radius $C^{-1}\log{r_j}$--- mapping to each point of $S^*M$ (see eg. \cite{LaxPhillips}), these contributions will be negligible with respect to the final $1/\log{r_j}$ bound.  Thus we must  estimate the projection of $\kappa_{\xi,j}\Big|_{U_j^{\text{up}}}$ to $V$.

Recall that the point
$$(x,t,\theta) = \begin{pmatrix} 1 & x\\ & 1 \end{pmatrix} \begin{pmatrix} e^{t/2} & 0\\ 0 & e^{-t/2} \end{pmatrix}. k_\theta = \begin{pmatrix} e^{t/2} & 0\\ 0 & e^{-t/2} \end{pmatrix} \begin{pmatrix} 1 & e^tx\\ & 1 \end{pmatrix} .k_\theta \in U_j^{\text{up}}$$
and so for a point given, in local $ANK$ coordinates around $p\in\xi$, by the matrix $\begin{pmatrix} e^{\tau/2} & \\ & e^{-\tau/2}\end{pmatrix}\begin{pmatrix} 1 & y\\ & 1\end{pmatrix} \begin{pmatrix} \cos\theta & -\sin\theta\\ \sin\theta & \cos\theta \end{pmatrix}\in V$ we have
\begin{eqnarray*}
\lefteqn{\bar{\kappa}_{\xi,j}\left(\begin{pmatrix} e^{\tau/2} & \\ & e^{-\tau/2}\end{pmatrix}\begin{pmatrix} 1 & y\\ & 1\end{pmatrix} \begin{pmatrix} \cos\theta & -\sin\theta\\ \sin\theta & \cos\theta \end{pmatrix}\right)}\\
 & = & \sum_{|k| l(\xi)\leq C^{-1}\log{r_j}} \kappa_{\xi,j}(ye^{kl(\xi)+\tau}, kl(\xi)+\tau,\theta) +  O(r_j^{-1/C})
\end{eqnarray*}
trivially estimating the contributions from $\kappa_{\xi,j}$ outside $U_j^{\text{up}}$ by $O(r_j^{-1/C})$,  as discussed above.  Therefore we  find by Proposition~\ref{spec integral} and the definition of $g(t)$ in section~\ref{kernel}
\begin{eqnarray*}
\lefteqn{ ||\bar\kappa_{\xi,j}||_{L^2(V)}^2\cdot \Big(1+ O(r_j^{-1/C})\Big)}\\
 & \gtrsim & \int_{\tau=0}^{l(\xi)} \int_{|y|\leq \eta} \int_{|\theta|\leq r^{-5/C}} \left(\sum_{|k| l(\xi)\leq C^{-1}\log{r_j}} \kappa_{\xi,j}(ye^{kl(\xi)+\tau}, kl(\xi)+\tau,\theta)	\right)^2 d\tau dy d\theta\\
 & \gtrsim &  \int_{\tau=0}^{l(\xi)} \left(\frac{\cos(r_j\tau)}{\log{r_j}}\right)^2
 \int_{y=0}^{\eta} \left(\sum_{|k| l(\xi)\leq \frac{1}{2}C^{-1}\log{r_j}} \chi(ye^{kl(\xi)+\tau})e^{\frac{1}{2}(kl(\xi)+\tau)}
 	\right)^2
 dyd\tau
 \end{eqnarray*}
since $\int_{|\theta|\leq r_j^{-5/C}} F_L^2(2\theta)d\theta \sim 1$ because $L\geq r_j^{10/C}$, the symmetry in $\chi$ allows us to integrate only over the positive values of $y$ at the expense of a constant factor of $2$, and restricting the range of $k$ decreases the overall value of the integrals, while allowing us to replace $\chi(\frac{C}{\log{r_j}}(kl(\xi)+\tau) / \cosh(\frac{C\pi(kl(\xi)+\tau)}{2\log{r_j}}) \gtrsim 1$ in the definition of $g(kl(\xi)+\tau)$on this range.  Note furthermore that since $r_j$ is a whole multiple of $2\pi l(\xi)$, we have $\cos(r_j(kl(\xi)+\tau)) = \cos(r_j\tau)$, and can be pulled out of the sum over $k$; here we use our assumption that $\{r_j\}$ is a sequence of  frequencies tailored to the geodesic $\xi$.

Moreover, further restricting the range of $y$ and noting that $\chi \equiv 1$ on $[-\frac{1}{2}, \frac{1}{2}]$, and that the $\tau$-dependence inside the sum can be abandoned at the expense of another constant depending only on $l(\xi)$,  we find
\begin{eqnarray*}
\lefteqn{ ||\bar\kappa_{\xi,j}||_{L^2(V)}^2\cdot \Big(1+ O(r_j^{-1/C})\Big)}\\
 & \gtrsim & \frac{1}{\log^2(r_j)}
\int_{y=0}^{\eta} \left(\sum_{|k| l(\xi) \leq \frac{1}{2}C^{-1}\log{r_j}} \chi(ye^{kl(\xi)})e^{\frac{1}{2}kl(\xi)}		\right)^2 dy\\
& \gtrsim & \frac{1}{\log^2(r_j)}\int_{y=r_j^{-1/2C}}^{\eta} \left(\sum_{-\frac{1}{2}C^{-1}\log{r_j} \leq kl(\xi) \leq -\log(2y)}e^{\frac{1}{2}kl(\xi)}		\right)^2 dy\\
&\gtrsim & \frac{1}{\log^2(r_j)}\int_{y=r_j^{-1/2C}}^{\eta} \left( e^{-\log(2y)/2}	\right)^2 dy\\
& \gtrsim & \frac{1}{\log^2(r_j)}\int_{y=r_j^{-1/2C}}^{\eta} \frac{1}{2y} dy\\
& \gtrsim &\frac{1}{\log^2(r_j)} (\log\eta + \frac{1}{2C}\log{r_j}) \gtrsim 1/C\log{r_j}
\end{eqnarray*}
where we have used the fact that $ye^{kl(\xi)}\leq \frac{1}{2}$ when $kl(\xi)\leq -\log(2y)$, and have restricted the integral to the range $y\geq r_j^{-1/2C}$ to ensure that $-\log(2y)\leq \frac{1}{2}C^{-1}\log{r_j}$ in changing the sum over $k$.

Since $||\bar{\kappa}_{\xi,j}||_{L^2(U)} \geq ||\bar{\kappa}_{\xi,j}||_{L^2(V)}$, the result follows.  $\Box$

\begin{Proposition}\label{L2 small}
We have for all $r_j$ sufficiently large
$$||\bar{\kappa}_{\xi, j}||_{L^2(S^*M)}^2 \lesssim 1/\log{r_j}$$
Combined with Proposition~\ref{Uj big}, this implies
$$||\bar{\kappa}_{\xi, j}||_{L^2(S^*M)}^2 \geq ||\bar{\kappa}_{\xi, j}||_{L^2(U)}^2 \gtrsim  1/\log{r_j} \gtrsim ||\bar{\kappa}_{\xi, j}||_{L^2(S^*M)}^2$$
for any neighborhood $U$ of $\xi$; thus any quantum limit $\mu$ of the probability measures 
$\mu_j ( f) = \frac{1}{||\bar{\kappa}_{\xi,j}||_2^2}\int_{S^*M} f |\bar{\kappa}_{\xi,j}|^2dVol$
satisfies $\mu(\xi)>0$. 
\end{Proposition}

{\em Proof:}  Again, since $\kappa_{\xi, j} \lesssim r^{-3/C}$ off of $U_j$ by Corollary~\ref{kappa small}, and the projection of $\text{supp}(\kappa_{\xi,j})$ to $S^*M$ is $O(r_j^{2/C})$-to-one, we may restrict our attention to the projection of $U_j$ to $S^*M$; i.e. 
$$||\bar{\kappa}_{\xi,j}||_{L^2(S^*M)} \leq ||\bar{\kappa}_{\xi,j}||_{L^2(\bar{U}_{j}^{\text{up}})} + ||\bar{\kappa}_{\xi,j}||_{L^2(\bar{U}_{j}^{\text{down}})} + O(r_j^{-1/C})$$
where $\bar{U}_j^{\text{up}}$ and $\bar{U}_j^{\text{down}}$ are the respective projections to $S^*M$ of $U_j^{\text{up}}$ and $U_j^{\text{down}}$, and we can evaluate the two norms on the right separately; since they are essentially identical, we will do the $\bar{U}_j^{\text{up}}$ term.  Note that
$$\bar{U}_j^{\text{up}} \subset \left\{ 	\xi.\begin{pmatrix}1 & y\\ & 1\end{pmatrix}	\begin{pmatrix} \cos\theta & -\sin\theta\\ \sin\theta & \cos\theta\end{pmatrix}   : |y|\leq r_j^{1/C}, \quad |\theta|\leq r_j^{-5/C} 	\right\}$$

We follow the calculations in the proof of Proposition~\ref{Uj big}; this time, though, we trivially estimate all $\chi$ terms, as well as the $\cos/\cosh$ term, and $\int_\theta F_L^2(2\theta)d\theta$, by $1$.  We recall however that $\chi$ is supported on $[-1,1]$, and so 
\begin{eqnarray*}
\lefteqn{ ||\bar\kappa_{\xi,j}||_{L^2(\bar{U}_j^{\text{up}})}^2}\\
 & \lesssim & \int_{\tau=0}^{l(\xi)} \int_{|y|\leq r_j^{1/C}} \int_{|\theta|\leq r_j^{-5/C}} \left(\sum_{|k| l(\xi)\leq C^{-1}\log{r_j}} \kappa_{\xi,j}(ye^{kl(\xi)+\tau}, kl(\xi)+\tau,\theta)	\right)^2 d\tau dy d\theta\\
 & \lesssim & \frac{1}{\log^2(r_j)} \int_{\tau=0}^{l(\xi)}
 \int_{y=0}^{r_j^{1/C}} \left(\sum_{|k| l(\xi)\leq C^{-1}\log{r_j}} \chi(ye^{kl(\xi)+\tau})e^{\frac{1}{2}(kl(\xi)+\tau)}
 	\right)^2
 dyd\tau\\
& \lesssim &  \frac{1}{\log^2(r_j)}
 \int_{y=r_j^{-1/C}}^{r_j^{1/C}} \left(\sum_{-C^{-1}\log{r_j}\leq kl(\xi)\leq -\log{y}} e^{\frac{1}{2}kl(\xi)}
 	\right)^2dy \\
	& & \quad + \quad \frac{1}{\log^2(r_j)}\int_{y=0}^{r_j^{-1/C}}  \left(\sum_{|k| l(\xi)\leq C^{-1}\log{r_j}} e^{\frac{1}{2}kl(\xi)}
 	\right)^2dy\\
 & \lesssim & \frac{1}{\log^2(r_j)}\left(\int_{y=r_j^{-1/C}}^{r_j^{1/C}} e^{-\log{y}} dy + \int_{y=0}^{r_j^{-1/C}} r_j^{1/C} dy\right)\\
 & \lesssim & \frac{1}{\log^2(r_j)}(C^{-1} \log{r_j} + 1) \lesssim 1/\log{r_j}
\end{eqnarray*}
as required.  The $\bar{U}_j^{\text{down}}$ term is identical, using the analogous trivial estimate of $\int_\theta F_L^2(\pi-2\theta)d\theta \leq 1$.
$\Box$

Finally, we must verify that $\bar{k}_{\xi,j}$ is indeed an $O\left(\frac{1}{\log{r_j}}\right)$-quasimode.  From the above calculations in Proposition~\ref{L2 small}, it is most convenient to do this via the estimates for ${\kappa}_{\xi,j}$.
\begin{Proposition}\label{quasimode}
We have
$$\left\|\left(\Delta + \frac{1}{4}+r_j^2	\right) \bar{k}_{\xi,j}\right\|_2 \lesssim \frac{r_j}{\log{r_j}}||\bar{k}_{\xi, j}||_2$$
and thus $\{\bar{k}_{\xi,j}\}$ is a sequence of $O\left(\frac{1}{\log{r_j}}\right)$-quasimodes.
\end{Proposition}

{\em Proof:}  Since $||\bar{k}_{\xi, j}||_{L^2(M)} = ||\bar{\kappa}_{\xi, j}||_{L^2(S^*M)}$, we can replace the right-hand side with $\frac{r_j}{\log{r_j}}||\bar{\kappa}_{\xi,j}||_{L^2(S^*M)}$.  For the left-hand side, we first compute on $\mathbb{H}$
$$\left(\Delta + \frac{1}{4}+r_j^2	\right)k_{\xi,j}  =  \int_s (r_j^2-s^2)k_{\xi,j}^{(s)} s h(s)\tanh(\pi s) ds$$
Since the operators $R^n$ and $R^{-n}$ defining the microlocal lift in Lemma~\ref{asymp to op} act in each irreducible representation of spectral parameter $s$, and $(\Delta + \frac{1}{4}+r_j^2)$ also commutes with $\Gamma$-translations on the left, we can  replace the left-hand side of Proposition~\ref{quasimode} with 
\begin{eqnarray*}
\left\|\left(\Delta + \frac{1}{4}+r_j^2	\right) \bar{k}_{\xi,j}\right\|_{L^2(M)} 
& = & \left\| \sum_{\gamma\in\Gamma} \int_s(r_j^2-s^2) k_{\xi,j}^{(s)}(\gamma\cdot)sh(s)\tanh{\pi s}ds\right\|_{L^2(M)}\\
& = & \left\| \sum_{\gamma\in\Gamma} \int_s(r_j^2-s^2) \kappa_{\xi,j}^{(s)}(\gamma\cdot)sh(s)\tanh{\pi s}ds\right\|_{L^2(S^*M)}
\end{eqnarray*}
and then we may perform all calculations in terms of the microlocal lifts $\kappa_{\xi,j}^{(s)}$, where we have already estimated the stationary phase asymptotics.

Indeed, by Proposition~\ref{stationary phase estimate}, and following the proof of Proposition~\ref{spec integral},
\begin{eqnarray*}
\lefteqn{\int_s (r_j^2-s^2) \kappa^{(s)}_{\xi,j} (x,t,\theta) sh(s)\tanh(\pi s)ds}\\
& = & \pi\chi(x)e^{t/2} \int_s (r_j^2-s^2)h(s)\Big(	e^{ist}F_L(2\theta) + e^{-ist}F_L(\pi-2\theta)	\Big) ds + O(r^{100/C-1})\\
& = & \pi\chi(x) e^{t/2} \Big(F_L(2\theta)  + F_L(\pi-2\theta)\Big) \cdot (r_j^2g(t) + g''(t)) + O(r^{100/C-1})
\end{eqnarray*}
and it remains to estimate the function $r_j^2g(t) + g''(t)$.

Recall that
\begin{eqnarray*}
g(t) & = & \frac{C}{4\log{r_j}}\cdot \frac{\cos(tr_j)}{\cosh{\left(\frac{C\pi t}{2\log{r_j}}\right)}} \chi\left(\frac{C}{\log{r_j}}t\right)\\
& = & \cos(tr_j) \cdot K H(Kt)
\end{eqnarray*}
where $K_j:=\frac{C}{2\log{r_j}}$, and $H(\xi) := \frac{1}{2}\chi(2\xi)/\cosh{\pi \xi}$.
Thus,
\begin{eqnarray*}
g''(t) + r_j^2g(t) & = & -r_j\sin(tr_j)\cdot K^2H'(Kt) + O(K^3H''(Kt))\\
& = & -(r_jK) \cdot \sin(tr_j) KH'(Kt) + O\left(\frac{1}{\log^3{r_j}}	\right)
\end{eqnarray*}
Now $\sin(tr_j)$ and $H'$ are uniformly bounded in $r_j$, so that the calculations of Proposition~\ref{L2 small}--- where $g(t)$ was estimated trivially in the support of $H$--- may be applied equally well to the function $\sin{(r_jt)}KH'(Kt)$ in place of $g(t)=\cos(r_jt)KH(Kt)$, and thus we arrive at
\begin{eqnarray*}
\left\|\left(\Delta + \frac{1}{4}+r_j^2	\right) \bar{k}_{\xi,j}\right\|_2^2 & \lesssim & (r_jK)^2 \cdot \frac{1}{\log{r_j}}\\
& \lesssim & \left(\frac{r_j^2}{\log^2{r_j}}\right) ||\bar{\kappa}_{\xi,j}||_{L^2(S^*M)}^2\\
& \lesssim & \left(\frac{r_j^2}{\log^2{r_j}}\right) ||\bar{k}_{\xi,j}||_{L^2(M)}^2
\end{eqnarray*}
as required, using the estimate $||\bar{\kappa}_{\xi,j}||_2^2 \asymp 1/\log{r_j}$ from Proposition~\ref{L2 small} and the fact that $||\bar{\kappa}_{\xi,j}||_{L^2(S^*M)} = ||\bar{k}_{\xi,j}||_{L^2(M)}$.   $\Box$

\vspace{.2in}

Combining Propositions~\ref{Uj big}, \ref{L2 small}, and \ref{quasimode}, we have now proved Theorem~\ref{loc quasi} for $\left(\frac{C'}{\log{r}}\right)$-quasimodes, with a sufficiently large constant $C' =O( C)$; in the next section we will describe how to build on this case to reduce the constant from a large $C'$ to arbitrarily small $\epsilon>0$.

\section{Finer-scale Quasimodes}\label{big to small}

We now discuss how to use these $\left(\frac{C'}{\log{r_j}}\right)$-quasimodes $\bar{k}_{\xi, j}$, that concentrate at least $\delta_1>0$ of their mass on the geodesic $\xi$, to construct $\left(\frac{\epsilon}{\log{r}}\right)$-quasimodes that concentrate positive mass on $\xi$, for any $\epsilon>0$.  Naturally, in agreement with Conjecture~\ref{hypconj}, our lower bound on the mass concentration will shrink as $\epsilon\to 0$.

The idea is to start with an arbitrary $o\left(\frac{1}{\log{r_j}}\right)$-quasimode $\psi_j$ near $r_j$--- for example, an eigenfunction with spectral parameter near $r_j$--- whose microlocal lift $\Psi_j$ we may assume does not concentrate mass on $\xi$, and perturb it by $\delta_2 \bar{k}_{\xi,j}$.  The  sum $\psi_j + \delta_2 \bar{k}_{\xi,j}$ will be an $\frac{\epsilon}{\log{r_j}}$-quasimode if $\delta_2$ is small enough, and $\Psi_j + \delta_2\bar{\kappa}_{\xi,j}$ will still concentrate  $\delta_2 \delta_1$ of its mass on $\xi$.

In order for this to work, $\psi_j$ must have approximate eigenvalue near $r_j$; precisely, within $K/\log{r_j}$, for some fixed $K$ independent of $\epsilon$.  By the estimates of B\'erard \cite{Berard} for the remainder term in Weyl's Law, we know that there exists a $K$ depending only on the surface $M$, such that an interval of width $K/\log{r}$ around $r$ must contain at least one eigenvalue.  So we may indeed take $\psi_j$ to be such an eigenfunction, with spectral parameter $\tilde{r}_j$ satisfying
$$0\leq \tilde{r}_j - r_j \leq K/\log{r_j}$$
and set $\delta_2  = \epsilon/(2K+C')$.  Naturally, we normalize $||\psi_j||_2 = ||\bar{k}_{\xi,j}||_2$.

If (hypothetically) the microlocal lifts $\Psi_j$ of these $\psi_j$ were to concentrate  $\geq \frac{1}{2}\delta_2\delta_1$ of its mass on $\xi$, then there is nothing to prove; these $\psi_j$ are eigenfunctions, and thus in particular are $\frac{\epsilon}{\log{r_j}}$-quasimodes for any $\epsilon>0$.  Otherwise, we define $\tilde{k}_j = \psi_j + \delta_2 \bar{k}_{\xi,j}$, whereby its microlocal lift is
$$\tilde{\kappa}_j = \Psi_j + \delta_2  \bar{\kappa}_{\xi,j}$$

Since by the triangle inequality, for any sufficiently small neighborhood $U\supset\xi$ in $S^*M$ and $j$ sufficiently large, we have
$$||\tilde{\kappa}_j||_{L^2(U)} \geq ||\delta_2 \bar\kappa_{\xi,j}||_{L^2(U)} - ||\Psi_j||_{L^2(U)} \geq \frac{1}{2} \delta_2 \delta_1 $$
we see that $\tilde{\kappa}_j$ concentrates a positive proportion of its mass on the geodesic $\xi$.  It remains to show that $\tilde{k}_j$ is indeed an $\frac{\epsilon}{\log{r_j}}$-quasimode.

But this is an immediate consequence of the fact that we diluted $\bar{k}_{\xi,j}$ by a sufficiently large factor.  Observe that since $\psi_j$ is an eigenfunction of spectral parameter $\tilde{r}_j$ we have
\begin{eqnarray*}
||(\Delta + (\frac{1}{4}+\tilde{r}_j^2))\tilde{k}_j||_2 & \leq &  ||(\Delta + (\frac{1}{4}+\tilde{r}_j^2))\psi_j||_2 + \delta_2 ||(\Delta + (\frac{1}{4}+\tilde{r}_j^2)) \bar{k}_{\xi,j}||_2\\
& \leq & \delta_2 ||(\Delta + (\frac{1}{4}+\tilde{r}_j^2)) \bar{k}_{\xi,j}||_2\\
& \leq & \delta_2 \Big(||(\Delta + (\frac{1}{4}+r_j^2)) \bar{k}_{\xi,j}||_2 + |r_j^2 - \tilde{r}_j^2|\cdot ||\bar{k}_{\xi,j}||_2\Big)\\
& \leq & \delta_2 \left(\frac{r_jC'}{\log{r_j}} + \frac{2\tilde{r}_jK}{\log{\tilde{r}_j}}\right)||\bar{k}_{\xi,j}||_2\\
& \leq & \tilde{r}_j \frac{\epsilon}{\log{\tilde{r}_j}} ||\psi||_2
\end{eqnarray*}
since $\delta_2=\epsilon/(2K+C')$, and $0\leq \tilde{r}_j-r_j\leq K/\log{r_j}$, completing the proof of Theorem~\ref{loc quasi}.  $\Box$

\def\cprime{$'$}

\end{document}